\documentclass[10pt]{amsart}
\usepackage{epsfig}
\usepackage{amsmath, amssymb, euscript,  enumerate}
\usepackage{bbm}
\usepackage{color}
\usepackage{verbatim}

\newcommand{\ap}{\alpha}             
\newcommand{\bt}{\beta}
\newcommand{\gm}{\gamma}             \newcommand{\Gm}{\Gamma}
             \newcommand{\Dt}{\Delta}
\newcommand{\vep}{\varepsilon}
\newcommand{\zt}{\zeta}

\newcommand{\ld}{\lambda}            
\newcommand{\sm}{\sigma}             
\newcommand{\vp}{\varphi}
\newcommand{\om}{\omega}             
\newcommand{\vr}{\varrho}            \newcommand{\iy}{\infty}
            
\newcommand{\f}{\frac}

\newcommand{\fa}{{\mathfrak a}}

\newcommand{\fe}{{\mathfrak e}}

\newcommand{\BN}{{\mathbb N}}

\newcommand{\BR}{{\mathbb R}}

\newcommand{\cC}{{\mathcal C}}

\newcommand{\cK}{{\mathcal K}}

\newcommand{\cO}{{\mathcal O}}
\newcommand{\cP}{{\mathcal P}}

\newcommand{\cS}{{\mathcal S}}

\newcommand{\la}{\langle}          \newcommand{\ra}{\rangle}
\newcommand{\s}{\setminus}         
            \newcommand{\e}{\eta}
\newcommand{\pa}{\partial}        
       
    \newcommand{\ds}{\displaystyle}

 \newcommand{\pf }{\noindent{\it Proof. }}

\newcommand{\aee }{\text{\rm a.e.}} 
  \newcommand{\pv }{\text{\rm p.v.}}

\newcommand{\rL }{{\text{\rm L}}}

\newcommand{\rX}{{\text{\rm X}}}  
\newcommand{\rS}{{\text{\rm S}}}

\newtheorem{thm}[subsection]{Theorem}
\newtheorem{lemma}[subsection]{Lemma}

\numberwithin{equation}{section}

\title[ the Berezin-Li-Yau type inequality for nonlocal operators]
{A refinement of the Berezin-Li-Yau type inequality for nonlocal
elliptic operators}
\author{ Yong-Cheol Kim }
\begin{document}
\begin{abstract} In this paper, we prove a refinement of the Berezin-Li-Yau type inequality for a wider class of
nonlocal elliptic operators including the fractional Laplacians
$-(-\Delta^{\sm/2})$ restricted to a bounded domain $D\subset\BR^n$
for $n\ge 2$ and $\sm\in (0,2]$, which is optimal when $\sigma=2$ in
view of Weyl's asymptotic formula. In addition, we describe the
Berezin-Li-Yau inequality for the Laplacian $\Delta$ as the limit
case of our result as $\sm\to 2^-$.
\end{abstract}

\address{$\bullet$ Yong-Cheol Kim : Department of Mathematics Education, Korea University, Seoul 136-701,
Korea }

\email{ychkim@korea.ac.kr}

\maketitle


\section{Introduction}

Let $\cK_{\sm}$ be the class of all positive symmetric kernels $K$
satisfying the uniformly ellipticity assumption
\begin{equation}K(y)=K(-y)\ge
\frac{\ld\,c_{n,\sm}}{|y|^{n+\sigma}},\,\,0<\sm<2,
\end{equation}
for all $y\in\BR^n\s\{0\}$ and $mK\in L^1(\BR^n)$, where
$m(y)=\min\{1,|y|^2\}$ and $c_{n,\sm}$ is the constant given by
\begin{equation}c_{n,\sm}=\biggl(\int_{\BR^n}\f{1-\cos(\xi_1)}{|\xi|^{n+\sm}}\,d\xi\biggr)^{-1}.
\end{equation}
Then we consider the corresponding nonlocal elliptic operator
$\rL_K$ given by
\begin{equation*}\rL_K
u(x)=\f{1}{2}\,\,\pv\int_{\BR^n}\mu(u,x,y)K(y)\,dy \end{equation*}
where $\mu(u,x,y)=u(x+y)+u(x-y)-2 u(x)$. In this paper, we consider
the following eigenvalue problem
\begin{equation}\begin{cases}-\rL_K u=\nu\,u &\text{ in $D$ }\\\qquad\, u=0
&\text{ in $\BR^n\s D$, }
\end{cases}\end{equation} where $\sm\in(0,2)$, $n\ge 2$, $K\in\cK_{\sm}$
and $D\subset\BR^n$ is an open bounded set.

Let $\rX$ be the normed linear space of all Lebesgue measurable
functions $v$ on $\BR^n$ with the norm
\begin{equation}\|v\|_{\rX}=\|v\|_{L^2(D)}+\biggl(\int_{\cC^n_D}|v(x)-v(y)|^2 K(x-y)\,dx\,dy\biggr)^{\f{1}{2}}<\iy
\end{equation}
where $\cC^n_H=\BR^{2n}\s(H^c\times H^c)$ for $H\subset\BR^n$. Set
$\rX_0=\{v\in\rX:v=0\,\,\aee\text{ in $\BR^n\s D$ }\}.$ Since
$C^2_0(D)\subset\rX_0$, we see that $\rX$ and $\rX_0$ are not empty.
By [9], there is a constant $c>1$ depending only on $n,\ld,\sm$ and
$D$ such that
\begin{equation*}\int_{\cC^n_D}|v(x)-v(y)|^2 K(x-y)\,dx\,dy\le\|v\|^2_{\rX}\le
c\int_{\cC^n_D}|v(x)-v(y)|^2 K(x-y)\,dx\,dy
\end{equation*}
for any $v\in\rX_0$; that is,
$\ds\|v\|_{\rX_0}:=\biggl(\int_{\cC^n_D}|v(x)-v(y)|^2
K(x-y)\,dx\,dy\biggr)^{1/2}$ is a norm on $\rX_0$ equivalent to
(1.4). Moreover it is known [9] that $(\rX_0,\|\cdot\|_{\rX_0})$ is
a Hilbert space with inner product
\begin{equation}\la
u,v\ra_{\rX_0}:=\iint_{\cC^n_D}(u(x)-u(y))(v(x)-v(y))K(x-y)\,dx\,dy.
\end{equation}
From simple computation, we note that $\la u,v\ra_{\rX_0}=-\la\rL_K
u,v\ra_{L^2(D)}$ for all $u,v\in\rX_0$.

More precisely, we study the weak formulation of the problem (1.3)
given by
\begin{equation}\begin{cases}\la u,v\ra_{\rX_0}=\nu\la u,v\ra_{L^2(D)},\,\forall v\in\rX_0,\\\qquad\,\,\,\,\, u\in\rX_0.
\end{cases}\end{equation}
Then it is well-known [10] that there is a sequence
$\{\nu^{\sm}_i(D)\}_{i\in\BN}$ of eigenvalues of (1.6) with
$0<\nu^{\sm}_1(D)\le\nu^{\sm}_2(D)\le\cdots\le\nu^{\sm}_i(D)\le\cdots$
and $\lim_{i\to\iy}\nu^{\sm}_i(D)=\iy$ such that the set
$\{e_i\}_{i\in\BN}$ of eigenfunctions $e_i$ corresponding to
$\nu^{\sm}_i(D)$ is an orthonormal basis of $L^2(D)$ and an
orthogonal basis of $\rX_0$. Moreover, it turns out that
$e_{i+1}\in\cP_{i+1}$ and
\begin{equation}\nu^{\sm}_1(D)=\| e_1\|^2_{\rX_0}\,\,\text{ and }\,\,\nu^{\sm}_{i+1}(D)=\| e_{i+1}\|^2_{\rX_0}
\end{equation}
for any $i\in\BN$, where $\cP_{i+1}=\{u\in\rX_0:\la
u,e_j\ra_{\rX_0}=0,\,\forall j=1,2,\cdots,i\}$.

Originally, Weyl's asymptotic formula [12] for the Dirichlet
eigenvalue problem of the Laplacian
\begin{equation}\begin{cases}-\Dt u=\mu\,u &\text{ in $D$ }\\\qquad\, u=0
&\text{ in $\pa D$ }
\end{cases}\end{equation}
asserts that
\begin{equation}\mu_k(D)\sim\f{4\pi^2}{(|D|\,|B_1|)^{\f{2}{n}}}
k^{\f{2}{n}}\,\,\text{ as $k\to\iy$, }
\end{equation} where $|D|$ and $|B_1|$ denote the volumes of $D$ and
the unit ball $B_1$ in $\BR^n$, respectively. The relevant study on
the eigenvalue problem for the Laplacian has been done along this
line by P\'olya [7] and Lieb [4]. P. Li and S. T. Yau [3] proved the
following lower bound on the averages on the finite sums of
eigenvalues
\begin{equation*}\f{1}{k}\sum_{j=1}^k\mu_k(D)\ge\f{4n\pi^2}{(n+2)(|D|\,|B_1|)^{\f{2}{n}}}k^{\f{2}{n}}
\end{equation*}
for any domain $D\subset\BR^n$, which is sharp in terms of (1.10).
P. Kr\"oger [2] obtained a upper bound for the sums of the
eigenvalues depending on geometric properties of $D$. A. Melas [5]
improved their lower bound by using the {\it moment of inertia} of
$D$. Using the method based on his argument, we obtain a lower bound
on the averages on the finite sums of eigenvalues $\nu^{\sm}_k$ of
the eigenvalue problem (1.6).

\begin{thm} Let $D\subset\BR^n$ be a bounded open
set and  $\sm\in(0,2)$. If $\,\{\nu^{\sm}_k(D)\}_{k\in\BN}$ be the
sequence of eigenvalues of the above eigenvalue problem $(1.6)$ for
the nonlocal elliptic operators $\rL_K$ with $K\in\cK_{\sm}$, then
we have the estimate
\begin{equation*}\f{1}{k}\sum_{j=1}^k\nu^{\sm}_k(D)
\ge\f{\ld\,n(2\pi)^{\sm}}{(n+\sm)(|B_1|\,|D|)^{\f{\sm}{n}}}k^{\f{\sm}{n}}
+\f{\ld\,\sm(2\pi)^{\sm-2}}{48(n+\sm)(|B_1|\,|D|)^{\f{\sm-2}{n}}}\f{|D|}{[D]}k^{\f{\sm-2}{n}}
\end{equation*} where $[D]=\int_D|x|^2\,dx$ is the moment
of inertia of $D$ with mass center $0\in\BR^n$.
\end{thm}

In particular, if $K_0(y)=c_{n,\sm}|y|^{-n-\sm}$ with $\sm\in(0,2)$,
then $\rL_{K_0}=-(-\Delta^{\sm/2})$ is the fractional Laplacian and
it is well-known [6] that
\begin{equation}\lim_{\sm\to
2^-}-(-\Delta^{\sm/2})u=\Delta u\,\,\text{ and }\,\,\lim_{\sm\to
0^+}(-\Delta^{\sm/2})u=u
\end{equation} for any function $u$
in the Schwartz space $\cS(\BR^n)$. Also, the constant $c_{n,\sm}$
satisfies the following property [6];
\begin{equation*}\lim_{\sm\to
2^-}\f{c_{n,\sm}}{\sm(2-\sm)}=\f{1}{|B_1|}\,\,\text{ and
}\,\,\lim_{\sm\to 0^+}\f{c_{n,\sm}}{\sm(2-\sm)}=\f{1}{2 n |B_1|}.
\end{equation*}

If we consider the nonlocal operator $\rL_{K_0}$ corresponding to
$K_0(y)=c_{n,\sm}|y|^{-n-\sm}$ with $\sm\in(0,2)$, our result makes
it possible to recover the result obtained by S. Yildirim Yolcu and
T. Yolcu [13] as follows.

\begin{thm} Let $D\subset\BR^n$ be a bounded open
set and  $\sm\in(0,2)$. If $\,\{\nu^{\sm}_{0k}(D)\}_{k\in\BN}$ be
the sequence of eigenvalues of the above eigenvalue problem $(1.6)$
for the fractional Laplacians $-(-\Delta^{\sm/2})$, then we have the
estimate
\begin{equation*}\f{1}{k}\sum_{j=1}^k\nu^{\sm}_{0k}(D)
\ge\f{n(2\pi)^{\sm}}{(n+\sm)(|B_1|\,|D|)^{\f{\sm}{n}}}k^{\f{\sm}{n}}
+\f{\sm(2\pi)^{\sm-2}}{48(n+\sm)(|B_1|\,|D|)^{\f{\sm-2}{n}}}\f{|D|}{[D]}k^{\f{\sm-2}{n}}
\end{equation*} where $[D]=\int_D|x|^2\,dx$ is the moment
of inertia of $D$ with mass center $0\in\BR^n$.
\end{thm}

As in (1.10), we could look on the Laplacian $\Delta$ as the limit
of the fractional Laplacian $-(-\Delta^{\sm/2})$ as $\sm\to 2^-$.
Then our result implies an improvement (see A. D. Melas [5]) of the
results proved by F. A. Berezin [1] and P. Li and S.-T. Yau [3] as
follows.

\begin{thm} Let $D\subset\BR^n$ be a bounded open
domain. If $\,\{\mu_k(D)\}_{k\in\BN}$ be the sequence of eigenvalues
of the above eigenvalue problem $(1.8)$ for the Laplacian $\Delta$,
then we have the estimate
\begin{equation*}\f{1}{k}\sum_{j=1}^k\mu_k(D)
\ge\f{n(2\pi)^2}{(n+2)(|B_1|\,|D|)^{\f{2}{n}}}k^{\f{2}{n}}
+\f{1}{24(n+2)}\f{|D|}{[D]}
\end{equation*} where $[D]=\int_D|x|^2\,dx$ is the moment
of inertia of $D$ with mass center $0\in\BR^n$.
\end{thm}

\section{Preliminaries}

First of all, we furnish several fundamental lemmas which are useful
in proving our main theorem. Our proof follows in part the argument
of Melas [5], Li and Yau [3], and S. Yildirim Yolcu and T. Yolcu
[13].

\begin{lemma} If $\phi:[0,\iy)\to[0,1]$ is a Lebesgue measurable function satisfying
$\int_0^{\iy}\phi(t)\,dt=1$ and  $0<\sm<2$, then there exists some
$\e>0$ such that
\begin{equation}\int_\e^{\e+1}t^n\,dt=\int_0^{\iy}t^n\phi(t)\,dt\,\,\,\text{
and
}\,\,\,\int_\e^{\e+1}t^{n+\sm}\,dt\le\int_0^{\iy}t^{n+\sm}\phi(t)\,dt.
\end{equation}
\end{lemma}

\pf First of all, we claim that $\int_0^{\iy}t^n\phi(t)\,dt<\iy$.
Indeed, if $\int_0^{\iy}t^{n+\sm}\phi(t)\,dt<\iy$, then we easily
obtain that $\int_0^{\iy}t^n\phi(t)\,dt<\iy$ because
$\int_0^{\iy}\phi(t)\,dt=1$. In case that
$\int_0^{\iy}t^{n+\sm}\phi(t)\,dt=\iy$, we can derive that
$\int_2^{\iy}t^{n+\sm}\phi(t)\,dt=\iy$. We note that the set
$H=\{t\in[2,\iy):\phi(t)>0\}$ must not be Lebesgue measure zero;
otherwise, it must be true that
$\int_2^{\iy}t^{n+\sm}\phi(t)\,dt=0$, which is a contradiction. So
we see that
$$\int_2^{\iy} t^n\phi(t)\,dt=\int_H t^n\phi(t)\,dt<\int_H t^{n+\sm}\phi(t)\,dt=\int_2^{\iy}t^{n+\sm}\phi(t)\,dt=\iy.$$
This implies that $\int_0^{\iy}t^n\phi(t)\,dt<\iy$.

Since $(t^n-1)(\phi(t)-\mathbbm{1}_{[0,1]}(t))\ge 0$ for any
$t\in[0,\iy)$, it follows from integrating the inequality on
$[0,\iy)$ that $\int_0^1 t^n\,dt\le\int_0^{\iy}t^n\phi(t)\,dt<\iy$.
We note that $g(s)=\int_s^{s+1} t^n\,dt$ is increasing on $[0,\iy)$
and $\lim_{s\to\iy}g(s)=\iy$. Thus there is an $\e>0$ such that
$g(\e)=\int_0^{\iy}t^n\phi(t)\,dt$.

We may also choose some $a,b\in(0,\iy)$ so that the function
\begin{equation*}h(t)=t^{n+\sm}-a t^n+b
\end{equation*}
satisfies $h(\e)=h(\e+1)=0$. Indeed, the equation $h'(t)=0$ has a
unique solution $(\f{an}{n+\sm})^{1/\sm}$ in $[0,\iy)$ at which  $h$
has the minimum value $-\f{a\sm}{n+\sm}(\f{an}{n+\sm})^{n/\sm}+b$.
Since $h$ is convex in $[0,\iy)$, we can select such $a,b\in(0,\iy)$
satisfying the condition $h(\e)=h(\e+1)=0$. Thus we conclude that
$h(t)<0$ for any $t\in(\e,\e+1)$ and $h(t)>0$ for any
$t\in[0,\iy)\s(\e,\e+1)$, and hence
$h(t)(\phi(t)-\mathbbm{1}_{(\e,\e+1)}(t))\ge 0$ for any
$t\in[0,\iy)$. Integrating this inequality on $[0,\iy)$, we easily
obtain the second result. \qed

\begin{lemma} The following inequality
\begin{equation}n t^{n+\sm}-(n+\sm)t^n s^{\sm}+\sm s^{n+\sm}\ge\sm
s^{n+\sm-2}(t-s)^2
\end{equation} always holds for any $s,t\in(0,\iy)$, $n\in\BN+1$ and
$\sm\in(0,2]$.
\end{lemma}

\pf If we set $\tau=t/s\in(0,\iy)$, then the inequality (2.2)
becomes
$$n\tau^{n+\sm}-(n+\sm)\tau^n+\sm\ge\sm(\tau-1)^2.$$
Consider the function
$p(\tau)=n\tau^{n+\sm}-(n+\sm)\tau^n+\sm-\sm(\tau-1)^2$. We write
$$p(\tau)=\tau\bigl(n\tau^{n+\sm-1}-(n+\sm)\tau^{n-1}-\sm\tau+2\sm\bigr):=\tau q(\tau).$$
Then we have that
$\,q'(\tau)=n(n+\sm-1)\tau^{n+\sm-2}-(n+\sm)(n-1)\tau^{n-2}-\sm$ and
\begin{equation*}\begin{split}q''(\tau)=\tau^{n-3}n(n+\sm-1)(n+\sm-2)\bigl(\tau^{\sm}-\tau_0^{\sm}\bigr)
\end{split}\end{equation*} where $\tau_0:=\bigl(\f{(n+\sm)(n-1)(n-2)}{n(n+\sm-1)(n+\sm-2)}\bigr)^{1/\sm}.$
The equation $q''(\tau)=0$ has a unique solution $\tau_0$ in
$(0,\iy)$ at which the function $q'(\tau)$ has the minimum value
$$q'(\tau_0)=\f{-\sm(n+\sm)(n-1)}{n+\sm-2}\biggl(\f{(n+\sm)(n-1)(n-2)}{n(n+\sm-1)(n+\sm-2)}\biggr)^{\f{n-2}{\sm}}-\sm<-\sm<0.$$
Since $\lim_{\tau\to 0^+}q'(\tau)=-\sm$ and $q'(1)=0$, we see that
the graph of $q'(\tau)$ is convex in $(0,\iy)$. Observing that
$\lim_{\tau\to 0^+}q(\tau)=2\sm$ and $q(1)=0$, this implies that the
graph of $q(\tau)$ is starting at the point $(0,2\sm)$ and going
down to the point $(1,0)$ convexly, and going up convexly right
after touching down to the point $(1,0)$. Hence we conclude that
$q(\tau)\ge 0$, and so $p(\tau)\ge 0$ for any $\tau\in(0,\iy)$. \qed

\begin{lemma} Let $n\in\BN+1$, $\vr,\bt\in(0,\iy)$ and $\sm\in(0,2]$. If
$\,\vp:[0,\iy)\to(0,\iy)$ is a decreasing absolutely continuous
function such that
\begin{equation}-\vr\le\vp'(t)\le 0\,\,\,\text{ and
}\,\,\,\int_0^{\iy}t^{n-1}\vp(t)\,dt:=\bt,
\end{equation}
then we have that
\begin{equation}\begin{split}\int_0^{\iy}t^{n+\sm-1}\vp(t)\,dt&\ge\f{1}{n+\sm}(n\bt)^{\f{n+\sm}{n}}\vp(0)^{-\f{\sm}{n}}\\
&\qquad\qquad\qquad+\f{\sm}{12
n(n+\sm)\vr^2}(n\bt)^{\f{n+\sm-2}{n}}\vp(0)^{\f{2n-\sm+2}{n}}.
\end{split}\end{equation}
\end{lemma}

\pf By considering the function $\vp(0)^{-1}\vp(\f{\vp(0)}{\vr}t)$,
we may assume that $\vr=1$ and $\vp(0)=1$. Without loss of
generality, we assume that
$\ap:=\int_0^{\iy}t^{n+\sm-1}\vp(t)\,dt<\iy$; otherwise, we have
already done. Set $\phi(t)=-\vp'(t)$ for $t\in[0,\iy)$. Then we have
that $0\le\phi(t)\le 1$ and $\int_0^{\iy}\phi(t)\,dt=\vp(0)=1$, and
moreover Lemma 2.1 and the integration by parts leads us to obtain
\begin{equation}\begin{split}\int_\e^{\e+1}t^{n+\sm}\,dt&\le\int_0^{\iy}t^{n+\sm}\phi(t)\,dt\\
&=\lim_{t\to\iy}\bigl(-t^{n+\sm}\vp(t)\bigr)+(n+\sm)\int_0^{\iy}t^{n+\sm-1}\vp(t)\,dt\\
&\le(n+\sm)\ap.
\end{split}\end{equation}
Applying the integration by parts again, by (2.5) we see that
\begin{equation}\begin{split}0\le\lim_{t\to\iy}t^{n+\sm}\vp(t)&=\bigl[t^{n+\sm}\vp(t)\bigr]^{\iy}_0\\
&=(n+\sm)\int_0^{\iy}t^{n+\sm-1}\vp(t)\,dt-\int_0^{\iy}t^{n+\sm}\phi(t)\,dt\\
&=(n+\ap)\ap-\int_0^{\iy}t^{n+\sm}\phi(t)\,dt<\iy.
\end{split}\end{equation}
Then we claim that $\lim_{t\to\iy}t^{n+\sm}\vp(t)=0\,$; indeed, if
$\gm:=\lim_{t\to\iy}t^{n+\sm}\vp(t)>0$, then given any
$\vep\in(0,\gm)$ there is some large $T>0$ such that
$\gm-\vep<t^{n+\sm}\vp(t)<\gm+\vep$ for all $t>T$, and thus we get
that
$$\iy=\int_T^{\iy}\f{\gm-\vep}{t}\,dt\le\int_0^{\iy}t^{n+\sm-1}\vp(t)\,dt<\iy,$$
which gives a contradiction. Hence, it follows from Lemma 2.1 and
the integration by parts that
\begin{equation}\int_\e^{\e+1}t^n\,dt=\int_0^{\iy}t^n\phi(t)\,dt=n\int_0^{\iy}t^{n-1}\vp(t)\,dt=n\bt.
\end{equation}
Integrating the inequality (2.2) on $[\e,\e+1]$, it follows from
(2.6) and (2.7) that
\begin{equation*}\begin{split}n(n+\sm)\ap-n(n+\sm)s^{\sm}\bt+\sm
s^{n+\sm}&\ge\sm s^{n+\sm-2}\int_\e^{\e+1}(t-s)^2\,dt\\
&\ge\sm s^{n+\sm-2}\int_{-1/2}^{1/2}t^2\,dt=\f{\sm}{2}\,s^{n+\sm-2}.
\end{split}\end{equation*}
Selecting $s=(n\bt)^{1/n}$, we obtain that
$$\ap\ge\f{1}{n+\sm}(n\bt)^{\f{n+\sm}{n}}+\f{\sm}{12
n(n+\sm)}(n\bt)^{\f{n+\sm-2}{n}}.$$ Therefore we complete the proof.
\qed

\section{Proof of Theorem 1.1}

In this section, we shall prove Theorem 1.1 by applying lemmas
obtained in the previous section.

Let $D\subset\BR^n$ be a bounded open domain and $D^*$ be its {\it
symmetric rearrangement} given by
$$D^*=\{x\in\BR^n:|x|<(|D|/|B_1|)^{1/n}\}.$$
That is, $D^*$ is the open ball with the same volume as $D$ and
center $0\in\BR^n$. Since $|x|^2$ is radial and increasing, the
moment of inertia of $D$ with mass center $0\in\BR^n$ has the lower
bound as follows;
\begin{equation}[D]=\int_D|x|^2\,dx\ge\int_{D^*}|x|^2\,dx=\f{n|D|}{n+2}\biggl(\f{|D|}{|B_1|}\biggr)^{\f{2}{n}}.
\end{equation}

Let $\{e_i\}_{i\in\BN}$ be the set of eigenfunctions $e_i$ of (1.7)
corresponding to eigenvalues $\nu^{\sm}_k(D)$ which is an
orthonormal basis of $L^2(D)$ and an orthogonal basis of $\rX_0$.
Then we consider the Fourier transform of each eigenfunction
$e_i(x)$ given by
\begin{equation*}\widehat e_i(\xi)=\f{1}{(2\pi)^{n/2}}\int_{\BR^n}e^{i\la x,\xi\ra}e_i(x)\,dx=\la \fe_\xi,e_i\ra_{L^2(D)}
\end{equation*}
where $\fe_\xi(x)=(2\pi)^{-n/2}e^{i\la x,\xi\ra}$. By Parseval's
formula and Plancherel theorem, we see that the set $\{\widehat
e_i\}_{i\in\BN}$ is orthonormal in $L^2(\BR^n)$. Since
$\{e_i\}_{i\in\BN}$ is an orthonormal basis of $L^2(D)$, it follows
from Bessel's inequality that
\begin{equation}\sum_{i=1}^k|\widehat
e_i(\xi)|^2\le\|\fe_{\xi}\|_{L^2(D)}=\f{|D|}{(2\pi)^n}
\end{equation}
for any $\xi\in\BR^n$ and $k\in\BN$. From standard analysis, we have
that
\begin{equation*}\nabla\widehat e_i(\xi)=\la ix\fe_{\xi},e_i\ra_{L^2(D)}
:=\bigl(\la ix_1\fe_{\xi},e_i\ra_{L^2(D)},\cdots,\la
ix_n\fe_{\xi},e_i\ra_{L^2(D)}\bigr).
\end{equation*}
Applying Bessel's inequality again, we obtain that
\begin{equation}\sum_{i=1}^k|\nabla\widehat
e_i(\xi)|^2\le\|ix\fe_{\xi}\|_{L^2(D)}=\f{[D]}{(2\pi)^n}
\end{equation}
for any $\xi\in\BR^n$ and $k\in\BN$. From (1.7) and Parseval's
formula, we have the estimate
\begin{equation}\begin{split}\nu^{\sm}_i(D)&=\|e_i\|^2_{\rX_0}=\la -\rL_K
e_i,e_i\ra_{L^2(D)}\\&=\la-\widehat{\rL_K e_i},\widehat
e_i\ra_{L^2(D)}=\int_{\BR^n}s(\xi)|\widehat e_i(\xi)|^2\,d\xi
\end{split}\end{equation}
where $s(\xi)=\int_{\BR^n}(1-\cos\la y,\xi\ra)K(y)\,dy$. Here we
note that $1-\cos\la y,\xi\ra\ge 0$. If we choose a matrix
$M\in\cO(n)$ such that $M e^1=\xi/|\xi|$ where
$e^1=(1,0,\cdots,0)\in\BR^n$, then by (1.1) we get the estimate
\begin{equation}\begin{split}s(\xi)&\ge\ld\,c_{n,\sm}\int_{\BR^n}\f{1-\cos\la
|\xi|y,\xi/|\xi|\ra}{|y|^{n+\sm}}\,dy\\
&=\ld\,c_{n,\sm}|\xi|^{\sm}\int_{\BR^n}\f{1-\cos\la \zt,M
e^1\ra}{|\zt|^{n+\sm}}\,d\xi\\
&=\ld\,c_{n,\sm}|\xi|^{\sm}\int_{\BR^n}\f{1-\cos\la \zt,
e^1\ra}{|\zt|^{n+\sm}}\,d\xi=\ld |\xi|^{\sm}.
\end{split}\end{equation}
If we set $G_k(\xi)=\sum_{i=1}^k|\widehat e_i(\xi)|^2$, then by
(3.2), (3.3) and Schwarz inequality, we have that $0\le
G_k(\xi)\le(2\pi)^{-n}|D|$. Also we observe that $x^{\fa}e_i\in
L^1(\BR^n)$ for any $i=1,\cdots,k$ and multi-index
$\fa=(a_1,\cdots,a_n)\in(\BN\cup\{0\}) ^n$, where
$x^{\fa}:=x_1^{a_1}\cdots x_n^{a_n}$. By Riemann-Lebesgue lemma and
standard analysis, we see that $\widehat e_i\in C^{\iy}_0(\BR^n)$
for any $i=1,\cdots,k$, and $G_k\in C^{\iy}_0(\BR^n)$. Thus we get
that
\begin{equation}|\nabla G_k(\xi)|\le 2\biggl(\,\sum_{i=1}^k|\widehat e_i(\xi)|^2\biggr)^{\f{1}{2}}
\biggl(\,\sum_{i=1}^k|\nabla\widehat
e_i(\xi)|^2\biggr)^{\f{1}{2}}\le\f{2\sqrt{|D|[D]}}{(2\pi)^n}
\end{equation}
for any $\xi\in\BR^n$, and moreover $\int_{\BR^n}G_k(\xi)\,d\xi=k$
by Plancherel theorem and
\begin{equation}\sum_{i=1}^k\nu^{\sm}_i(D)\ge \ld
\int_{\BR^n}|\xi|^{\sm}G_k(\xi)\,d\xi
\end{equation}
by (3.4) and (3.5).  Let $G_k^*(\xi)=\vp(|\xi|)$ be the symmetric
decreasing arrangement of $G_k$. Then it follows from Lemma 1.E. in
[11] that $\vp$ is absolutely continuous in $[0,\iy)$. For $\tau\ge
0$, we set
$\om(\tau)=|\{\xi\in\BR^n:G_k^*(\xi)>\tau\}|=|\{\xi\in\BR^n:G_k(\xi)>\tau\}|$.

\begin{lemma} If $\,\vp$ is differentiable at $t\in(0,\iy)$, then
$\om$ is differentiable at $\vp(t)$ and moreover
$\om'(\vp(t))\vp'(t)=n|B_1|t^{n-1}$.
\end{lemma}

\pf Take any $t\in(0,\iy)$ at which $\vp$ is differentiable. Then we
have two possible cases; (i) there is an open interval
$I\subset(0,\iy)$ such that $t\in I$ and $\vp'=0$ in $I$, and (ii)
there is an interval $I\subset(0,\iy)$ such that $t\in I$ and
$\vp'<0$ in $I$.

In case of (i), it is easy to check that $\om'(\vp(t))=0$. In case
of (ii), by the property of the distribution function, we see that
$\om$ is continuous at $\vp(t)$. We note that
$\om(\vp(t))=|B_1|t^n$. Write $\Delta s=\vp(t+\Delta t)-\vp(t)$.
Then we have that
\begin{equation*}\begin{split}\f{\om(\vp(t)+\Delta s)-\om(\vp(t))}{\Delta
s}\,\f{\Delta s}{\Delta t}&=\f{\om(\vp(t+\Delta
t))-\om(\vp(t))}{\Delta t}\\
&=|B_1|\sum_{i=0}^{n-1}(t+\Delta t)^{n-1-i}(\Delta t)^i.
\end{split}\end{equation*}
Taking the limit in the above because $\vp$ is continuous at $t$,
this implies the required result. \qed

We continue the proof of Theorem 1.1. As in the above, there is
nothing to prove it, because $\vp'=0$ in $I$ in case of (i). So,
without loss of generality, we may assume that we are now in the
case (ii). By (3.1) and (3.7), we have that
\begin{equation}k=\int_{\BR^n}G_k(\xi)\,d\xi=\int_{\BR^n}G_k^*(\xi)\,d\xi=n|B_1|\int_0^{\iy}t^{n-1}\vp(t)\,dt
\end{equation}
and
\begin{equation}\begin{split}\sum_{i=1}^k\nu^{\sm}_i(D)&\ge\ld
\int_{\BR^n}|\xi|^{\sm}G_k(\xi)\,d\xi\\&\ge\ld
\int_{\BR^n}|\xi|^{\sm}G_k^*(\xi)\,d\xi=\ld
n|B_1|\int_0^{\iy}t^{n+\sm-1}\vp(t)\,dt. \end{split}\end{equation}

Since $\vp:[0,\iy)\to [0,(2\pi)^{-n}|D|]$ is decreasing by (3.2), it
follows from the coarea formula that
\begin{equation}\om(\tau)=\int_\tau^{(2\pi)^{-n}|D|}\int_{\rS_s}\f{1}{|\nabla
G_k(\xi)|}\,d\sm_s(\xi)\,ds
\end{equation}
where $\rS_s=\{\xi\in\BR^n:G_k(\xi)=s\}$ and $d\sm_s$ is the surface
measure on $\rS_s$. Thus by (3.10) and Lemma 3.1, we obtain that
\begin{equation}\begin{split}-n|B_1|t^{n-1}&=-\om'(\vp(t))\vp'(t)=\biggl(\int_{\rS_{\vp(t)}}\f{1}{|\nabla
G_k(\xi)|}\,d\sm_{\vp(t)}(\xi)\biggr)\vp'(t)\\
&\le\f{1}{\vr}\,\sm_{\vp(t)}(\rS_{\vp(t)})\,\vp'(t)=\f{1}{\vr}\,n|B_1|t^{n-1}\,\vp'(t)\le
0,
\end{split}\end{equation}
where $\vr=2(2\pi)^{-n}\sqrt{|D|[D]}$. Thus this implies that
$-\vr\le\vp'(t)\le 0$ for any $t\ge 0$. If we set $\bt=k/(n|B_1|)$
in (3.8), then by (3.9) and Lemma 2.3 we have that
\begin{equation}\begin{split}\f{1}{k}\sum_{i=1}^k\nu^{\sm}_i(D)&\ge\f{\ld\,n
}{n+\sm}\biggl(\f{k}{|B_1|}\biggr)^{\f{\sm}{n}}\vp(0)^{-\f{\sm}{n}}\\
&\qquad\qquad\qquad+\f{\ld\,\sm
}{12(n+\sm)\vr^2}\biggl(\f{k}{|B_1|}\biggr)^{\f{\sm-2}{n}}\vp(0)^{\f{2n-\sm+2}{n}}.
\end{split}\end{equation}
For $t\in[0,(2\pi)^{-n}|D|]$, we set
$$h(t)=\f{\ld \,n}{n+\sm}\biggl(\f{k}{|B_1|}\biggr)^{\f{\sm}{n}}t^{-\f{\sm}{n}}
+\f{\ld\,\sm}{12(n+\sm)\vr^2}\biggl(\f{k}{|B_1|}\biggr)^{\f{\sm-2}{n}}t^{\f{2n-\sm+2}{n}}.$$
Differentiating $h(t)$ once, we get that
$$h'(t)=\f{\ld\,\sm}{n+\sm}\biggl(\f{k}{|B_1|}\biggr)^{\f{\sm}{n}}t^{-\f{\sm}{n}-1}g(t)$$
where
$g(t)=-1+\f{2n-\sm+2}{12n(n+\sm)\vr^2}(\f{k}{|B_1|})^{-2/n}t^{\f{n-\sm+2}{n}}$.
Since $g$ is increasing on $[0,(2\pi)^{-n}|D|]$, we obtain that
$$g(t)\le g((2\pi)^{-n}|D|)=-1+\f{(n+2)(2n-\sm+2)|B_1|^{4/n}}{192 n^2\pi^2
k^{2/n}}\le-1+\f{20}{192}\le 0$$ from the fact that $\sm\in(0,2]$,
$|B_1|=\f{2\pi^{n/2}}{n\Gm(\f{n}{2})}$ and
$\Gm(\f{n}{2})\ge\Gm(\f{1}{2})=\sqrt{\pi}$ for all $n\in\BN$. Thus
$h$ is decreasing on $[0,(2\pi)^{-n}|D|]$, a lower bound in (3.12)
can be obtained by replacing $\vp(0)$ by $(2\pi)^{-n}|D|$ as
follows;
\begin{equation*}\begin{split}\f{1}{k}\sum_{i=1}^k\nu^{\sm}_i(D)&\ge
h\bigl((2\pi)^{-n}|D|\bigr)\\
&=\f{\ld\,n(2\pi)^{\sm}}{(n+\sm)(|B_1|\,|D|)^{\f{\sm}{n}}}k^{\f{\sm}{n}}
+\f{\ld\,\sm(2\pi)^{\sm-2}}{48(n+\sm)(|B_1|\,|D|)^{\f{\sm-2}{n}}}\f{|D|}{[D]}k^{\f{\sm-2}{n}}.
\end{split}\end{equation*}
Therefore we complete the proof. \qed


\end{document}